
\tolerance=10000
\raggedbottom

\baselineskip=15pt
\parskip=1\jot

\def\sk{\vskip 3\jot}

\def\heading#1{\vskip3\jot{\noindent\bf #1}}
\def\label#1{{\noindent\it #1}}


\def\ref#1;#2;#3;#4;#5.{\item{[#1]} #2,#3,{\it #4},#5.}
\def\refinbook#1;#2;#3;#4;#5;#6.{\item{[#1]} #2, #3, #4, {\it #5},#6.} 
\def\refbook#1;#2;#3;#4.{\item{[#1]} #2,{\it #3},#4.}


\def\({\bigl(}
\def\){\bigr)}

\def\la{\lambda}

\def\De{\Delta}

\def\Ex{{\rm Ex}}
\def\Var{{\rm Var}}

{
\pageno=0
\nopagenumbers
\rightline{\tt egp.iii.arxiv.tex}
\vskip1in

\centerline{\bf The $M/M/\infty$ Service System with Ranked Servers in Heavy Traffic}
\vskip0.5in

\centerline{Patrick Eschenfeldt}
\centerline{\tt peschenfeldt@hmc.edu}
\sk

\centerline{Ben Gross}
\centerline{\tt bgross@hmc.edu}
\sk

\centerline{Nicholas Pippenger}
\centerline{\tt njp@math.hmc.edu}
\sk

\centerline{Department of Mathematics}
\centerline{Harvey Mudd College}
\centerline{1250 Dartmouth Avenue}
\centerline{Claremont, CA 91711}
\vskip0.5in

\noindent{\bf Abstract:}
We consider an $M/M/\infty$ service system in which an arriving customer is served by the 
first idle server in an infinite sequence $S_1, S_2, \ldots$ of servers.
We determine the first two terms in the asymptotic expansions of the moments of $L$ as $\la\to\infty$, where $L$ is the index of the server $S_L$ serving a newly arriving customer in equilibrium, and $\la$ is the ratio of the arrival rate to the service rate.
The leading terms of the moments show that $L/\la$ tends to a uniform distribution on $[0,1]$.
\vskip0.5in
\leftline{{\bf Keywords:} Queueing theory,  asymptotic expansions.}
\sk
\leftline{{\bf Subject Classification:} 60K26, 90B22}

\vfill\eject
}

\heading{1.  Introduction}

We consider a stream of customers, with independent exponentially distributed interarrival times,
arriving at rate $\la$ to an infinite sequence $S_1, S_2, \ldots$ of servers.
Each arriving customer engages the server $S_l$ having the lowest index among currently idle servers, and renders that server busy for an independent exponentially distributed service time with mean $1$.
This stochastic service system, which is conventionally denoted $M/M/\infty$, has been extensively studied in the limit $\la\to\infty$; see Newell [N].
We shall be interested in a question mentioned only tangentially by Newell:
what is the distribution of
the random variable $L$ defined as the index of the server
$S_L$ serving a newly arriving customer when the system is in equilibrium?
Newell [N, p.~9] states that $L$ ``is approximately uniformly distributed over the interval'' $[1, \la]$,
basing this assertion on the approximation
$$\Pr[L>l] \approx 
\cases{
1 - \displaystyle{l\over\la}, &if $l<\la$, \cr
0, &if $l>\la$. \cr
} \eqno(1.1)$$
But no error bounds are given for this or other approximations stated by Newell, and not even 
the fact that the first moment has the asymptotic behavior
$$\Ex[L]\sim{{\la\over 2}} \eqno(1.2)$$
that it would have under the uniform distribution is established rigorously.
Our goal in this paper is to give a rigorous version of (1.1) that will suffice to establish
not only (1.2), but also the next term,
$$\Ex[L] = {\la\over 2} + {1\over 2}\log \la + O(1), \eqno(1.3)$$
and more generally
$$\Ex[L^m] = {\la^m \over m+1} + {m\,\la^{m-1}\log\la\over 2} + O\left(\la^{m-1}\right) \eqno(1.4)$$
for $m\ge 1$.
In particular,
we have
$$\eqalign{
\Var[L] 
&= \Ex[L^2] - \Ex[L]^2 \cr
& = {\la^2\over 12} + {\la\log\la\over 2} + O(\la). \cr
}$$
Since the interval $[0,1]$ is bounded, formula (1.4) shows that the $m$-th moment of $L/\la$ tends to $1/(m+1)$ as $\la\to\infty$ for all $m\ge 1$, and thus suffices to show that the distribution of $L/\la$ tends to the uniform distribution on the interval  $[0,1]$.
We note that a problem that is in a sense dual to ours (finding the largest index of a busy server, rather than the smallest index of an idle server) has been treated by Coffman, Kadota and Shepp [C].

The key to our results is the probability $\Pr[L>l]$, which is simply the probability that the first $l$ servers $S_1,\ldots, S_l$ are all busy.
It is well known that this probability is given by the Erlang loss formula
$$\eqalign{
\Pr[L>l] 
&= {\la^l/l! \over \sum_{0\le k\le l} \la^k/k!} \cr
&= {1\over D_l}, \cr
}$$
where
$$D_l = \sum_{0\le k\le l} {l! \over (l-k)! \, \la^k} \eqno(1.5)$$
(see for example Newell [N, p.~3]).
The sum $D_l$ can be expressed as an integral,
$$D_l = \int_0^\infty \left(1 + {x\over\la}\right)^l \, e^{-x} \, dx$$
(see for example Newell [N, p.~7]), and most of Newell's analysis is based on such a representation. But we shall work directly with the expression of $D_l$ as the sum in (1.5).

We shall divide the range of summation in (1.5) into two parts.
The first, which we shall call the ``body'' of the distribution, will be $0\le k\le l_0 = \la-s$,
where $s = \sqrt{\la}$.
The second, which we shall call the ``tail'', will be $l > l_0$.
In Section 2, we shall derive an estimate for $\Pr[L>l]$ in the body, and in Section 3, we shall derive an estimate for the tail.
In Section 4, we shall combine these estimates to establish (1.4).
\sk

\heading{2. The Body}

In this section we shall establish the estimate
$$\Pr[L>l] = (1 - l/\la) +  {1 \over \la(1 - l/\la)}  + O\left({1\over\la}\right)
+ O\left({1\over \la^2(1-l/\la)^3}\right) \eqno(2.1)$$
for $l\le l_0 = \la - s$, where $s = \sqrt{\la}$.
We begin by using the principle of inclusion-exclusion to derive bounds on the denominator
$D_l$.

We begin with a lower bound.
Since
$$\eqalign{
l(l-1)\cdots(l-k+1) 
&\ge l^k  - \left(\sum_{0\le j\le k-1} j\right) l^{k-1} \cr
&= l^k - {k\choose 2}l^{k-1}, \cr
}$$
we have
$$\eqalign{
D_l
&= \sum_{0\le k\le l} {l(l-1)\cdots(l-k+1)\over\la^k} \cr
&\ge \sum_{0\le k\le l} \left({l\over\la}\right)^k  
- {1\over\la}\sum_{0\le k\le l} {k\choose 2}\left({l\over\la}\right)^{k-1}. \cr
}$$
For the first sum we have
$$\sum_{0\le k\le l} \left({l\over\la}\right)^k   = {1 + O\((l/\la)^l\)\over 1 - l/\la}.$$
We note that the logarithm of $(l/\la)^l$ has a non-negative second derivative for $l\ge 1$.
Thus $(l/\la)^l$ assumes its maximum in the interval $0\le l\le l_0$ for $l=0$,  $l=1$ or
$l=l_0$.
Its values there are $0$, $1/\la$ and $(1-s/\la)^{\la-s} = 
(1 - 1/\sqrt{\la})^{\la - \sqrt{\la}} \le e^{-\sqrt{\la} + 1}$, respectively.
As $\la\to\infty$, the largest of these values is $1/\la$, so we have 
$O\((l/\la)^l\) = O(1/\la)$ for $0\le l\le l_0$.
Thus the first sum is
$$\sum_{0\le k\le l} \left({l\over\la}\right)^k   = {1 + O(1/\la)\over 1 - l/\la}.$$
For the second sum we have
$$\sum_{0\le k\le l}{k\choose 2} \left({l\over\la}\right)^{k-1}   = {1 + O\(l^2(l/\la)^l\)\over (1 - l/\la)^3}.$$
The logarithm of $l^2(l/\la)^l$ has a non-negative second derivative for $l\ge 3$, so an argument similar to that used for the first sum shows that $O\(l^2(l/\la)^l\) = O(1/\la)$ for $0\le l\le l_0$.
Thus we have
$$\sum_{0\le k\le l}{k\choose 2} \left({l\over\la}\right)^{k-1}   = {1 + O(1/\la)\over (1 - l/\la)^3}$$
and the lower bound
$$D_l \ge {1 + O(1/\la)\over 1 - l/\la} - {1 + O(1/\la)\over \la(1 - l/\la)^3}. \eqno(2.2)$$

For an upper bound, we have
$$\eqalign{
l(l-1)\cdots(l-k+1) 
&\le l^k - \left(\sum_{0\le j\le k-1} j\right) l^{k-1} 
+ \left(\sum_{0\le i<j\le k-1} ij\right) l^{k-2} \cr
&\le l^k - {k\choose 2}l^{k-1} + {1\over 2}{k\choose 2}^2 l^{k-2} \cr
}$$
(because $\sum_{0\le i<j\le k-1} ij = \left(\left(\sum_{0\le j\le k-1} j\right)^2 - \sum_{0\le j\le k-1} j^2\right)\bigg/2
\le \left(\sum_{0\le j\le k-1} j\right)^2\big/2 = {k\choose 2}^2/2$).
Thus we have
$$D_l \le \sum_{0\le k\le l} \left({l\over\la}\right)^k  
- {1\over\la}\sum_{0\le k\le l} {k\choose 2}\left({l\over\la}\right)^{k-1}
+ {1\over2\la^2}\sum_{0\le k\le l} {k\choose 2}^2 \left({l\over\la}\right)^{k-2}.$$
For the third sum we have
$$\eqalign{
\sum_{0\le k\le l} {k\choose 2}^2 \left({l\over\la}\right)^{k-2}
&\le \sum_{k\ge 0} {k\choose 2}^2 \left({l\over\la}\right)^{k-2} \cr
&= O\left({1\over (1-l/\la)^5}\right). \cr
}$$
and thus the upper bound
$$D_l \le {1 + O(1/\la)\over 1 - l/\la} - {1 + O(1/\la)\over \la(1 - l/\la)^3}
 + O\left({1\over \la^2(1-l/\la)^5}\right).$$
 Combining this upper bound with the lower bound (2.2) yields
 $$D_l = {1 + O(1/\la)\over 1 - l/\la} - {1 + O(1/\la)\over \la(1 - l/\la)^3}
 + O\left({1\over \la^2(1-l/\la)^5}\right).$$
 
 To obtain $\Pr[L>l]$, we take the reciprocal of $D_l$:
 $$\eqalign{
 \Pr[L>l]
 &= \left({1 + O(1/\la)\over 1 - l/\la} - {1 + O(1/\la)\over \la(1 - l/\la)^3}
 + O\left({1\over \la^2(1-l/\la)^5}\right)\right)^{-1} \cr
  &= \(1 + O(1/\la)\)\;(1 - l/\la) \, \left(1 - {1 \over \la(1 - l/\la)^2}
 + O\left({1\over \la^2(1-l/\la)^4}\right)\right)^{-1} \cr
  &= \(1 + O(1/\la)\)\;(1 - l/\la) \, \left(1 + {1 \over \la(1 - l/\la)^2}
 + O\left({1\over \la^2(1-l/\la)^4}\right)\right) \cr
 &= \(1 + O(1/\la)\)\; \left((1 - l/\la) +  {1 \over \la(1 - l/\la)}  
 + O\left({1\over \la^2(1-l/\la)^3}\right)\right). \cr
 }$$
Observing that $O(1/\la)\,(1-l/\la) = O(1/\la)$ and
$O(1/\la)/\la(1-l/\la) = O\(1/\la^2(1-l/\la)^3\)$, we obtain (2.1).
\sk

\heading{3. The Tail}

In this section we shall establish the estimate
$$\Pr[L>l] = O(e^{-\la} \, \la^l / l!) \eqno(3.1)$$
for $l\ge \la - s$, where $s = \sqrt{\la}$.
To obtain an upper bound on $\Pr[L>l]$, we obtain a lower bound on $D_l$.
We have
$$\eqalignno{
D_l
&= \sum_{0\le k\le l} {l! \over (l-k)! \, \la^k} \cr
&\ge {l! \over \lfloor \la-s\rfloor! \, \la^{l-\lfloor \la-s\rfloor}} + \cdots
+ {l! \over \lfloor \la-2s\rfloor! \, \la^{l-\lfloor \la-2s\rfloor}}, &(3.2)\cr
}$$
because $l-\lfloor \la-s\rfloor \ge l- (\la-s) \ge 0$ by assumption and 
$\lfloor \la-2s\rfloor \ge 0$ for all sufficiently large $\la$.
There are $\lfloor \la-2s\rfloor - \lfloor \la-2s\rfloor + 1 \ge s$ terms
in the sum (3.2).
Furthermore, the smallest of these terms is the last, because its denominator contains factors of 
$\la$ where the preceding terms contain factors smaller than $\la$.
Thus we have
$$D_l \ge {s \, l! \over \lfloor \la-2s\rfloor! \, \la^{l-\lfloor \la-2s\rfloor}}.$$
For the factorial in the denominator of this bound, we shall use the estimate
$n! \le e \, \sqrt{n} \, e^{-n} \, n^n$, which holds for all $n\ge 1$
(because the trapezoidal rule underestimates the integral $\int_1^n \log x\,dx$ of the concave function $\log x$).
This estimate yields
$$D_l \ge {s \, l! \, e^{\lfloor \la-2s\rfloor} \over 
e\, \sqrt{\lfloor \la-2s\rfloor} \, \lfloor \la-2s\rfloor^{\lfloor \la-2s\rfloor} 
\, \la^{l-\lfloor \la-2s\rfloor}}. \eqno(3.3)$$
We have
$$e^{\lfloor \la-2s\rfloor} \ge e^{\la-2s-1},$$
$$\eqalign{
\lfloor \la-2s\rfloor^{\lfloor \la-2s\rfloor} 
&\le (\la-2s)^{\lfloor \la-2s\rfloor} \cr
&= \la^{\lfloor \la-2s\rfloor} \, (1-2s/\la)^{\lfloor \la-2s\rfloor} \cr
&\le \la^{\lfloor \la-2s\rfloor} \, (1-2s/\la)^{\la-2s-1} \cr
&\le \la^{\lfloor \la-2s\rfloor} \, e^{(-2s/\la)(\la-2s-1)} \cr
&\le \la^{\lfloor \la-2s\rfloor} \, e^{-2s+4s^2/\la + 1} \cr
&\le \la^{\lfloor \la-2s\rfloor} \, e^{-2s+5} \cr
}$$
and
$$\sqrt{\lfloor \la-2s\rfloor} \le s.$$
Substituting these bounds into (3.3) yields
$$D_l \ge {l! \, e^\la \over e^7 \, \la^l}.$$
Taking the reciprocal of this bound yields (3.1).
\vfill\eject

\heading{4. The Moments}

In this section we shall use (2.1) and (3.1) to prove (1.4).
We write
$$\eqalign{
\De_m(l) 
&= l^m - (l-1)^m \cr
&= m\,l^{m-1} + O(l^{m-2}) \cr
}$$
for the backward differences of the $m$-th powers of $l$.
Then partial summation yields
$$\eqalignno{
\Ex[L^m]
&= \sum_{l\ge 0} l^m \, \Pr[L=l] \cr
&= \sum_{l\ge 0} \De_m(l) \, \Pr[L>l] \cr
&= \sum_{l\ge 0} m\, l^{m-1} \, \Pr[L>l] + O\left(\sum_{l\ge 0} l^{m-2} \, \Pr[L>l]\right) &(4.1)\cr
}$$
This formula shows that we should evaluate sums of the form
$$T_n = \sum_{l\ge 0} l^n \, \Pr[L>l]. \eqno(4.2)$$
We shall show that
$$T_n = {\la^{n+1} \over (n+1)(n+2)} + {\la^n \, \log \la \over 2} + O(\la^n). \eqno(4.3)$$
Substitution of this formula into (4.1) will then yield (1.4).

We shall break the range of summation in (4.2) at $l_0 = \la - s$, where $s=\sqrt{\la}$,
using (2.1) for $0\le l\le l_0$ and (3.1) for $l>l_0$.
Summing the first term in (2.1),
we have
$$\eqalignno{
\sum_{0\le l\le l_0} l^n (1 - l/\la)
&= {1\over\la} \sum_{0\le l\le l_0} (\la\,l^n - l^{n+1}) \cr
&= {1\over\la} \left(\left({\la \, l_0^{n+1}\over n+1} + O(l_0^n)\right)
- \left({\la^{n+2} \over n+2} + O(l_0^{n+1})\right)\right) \cr
&= {1\over\la} \left(\left({\la \, (\la^{n+1} - (n+1)\la^n s)\over n+1} + O(\la^n)\right)
- \left({\la^{n+2}  - (n+2)\la^{n+1} s\over n+2} + O(\la^{n+1})\right)\right) \cr
&= {\la^{n+1} \over (n+1)(n+2)} + O(\la^n).  \cr
}$$
Summing the second term in (2.1), we have
$$\eqalignno{
\sum_{0\le l\le l_0} {l^n \over \la-l}
&= \sum_{s\le k\le\la} {(\la - k)^n \over k} \cr
&= \sum_{s\le k\le\la} \left({\la^n \over k} + O(\la^{n-1})\right) \cr
&= \la^n \log {\la\over s} + O(\la^n) \cr
&= {\la^n \log\la \over 2} + O(\la^n), \cr
}$$
where we have used $\sum_{1\le k\le n} 1/k = \log n + O(1)$.
Summing the third term in (2.1) of course yields $O(\la^n)$.
Summing the last term in (2.1), we have
$$\eqalignno{
\la \sum_{0\le l\le l_0} {l^n \over (\la - l)^3}
& = \la \sum_{s\le k\le\la} {(\la-k)^n \over k^3} \cr
& \le \la^{n+1} \sum_{s\le k\le\la} {1 \over k^3} \cr
& \le \la^{n+1} \sum_{k\ge s} {1 \over k^3} \cr
& = \la^{n+1} \left({2\over s^2} + O\left({1\over s^3}\right)\right) \cr
&= O(\la^n), \cr
}$$
where we have used $\sum_{k\ge n} 1/k^3 = 2/n^2 + O(1/n^3)$.
Combining these estimates, we obtain
$$\sum_{0\le l\le l_0} l^n \, \Pr[L>l] 
= {\la^{n+1} \over (n+1)(n+2)} + {\la^n \, \log \la \over 2} + O(\la^n). \eqno(4.4)$$
Finally, summing (3.1) we have
$$\eqalign{
\sum_{l>l_0} {l^n \, e^{-\la} \, \la^l \over l!} 
&\le \sum_{l\ge 0} {l^n \, e^{-\la} \, \la^l \over l!}  \cr
&= O(\la^n), \cr
}$$
because the summation on the right-hand side is the $n$-th moment of a Poisson random variable with mean $\la$, which is a polynomial of degree $n$ in $\la$.
Thus
$$\sum_{l > l_0} l^n \, \Pr[L>l] = O(\la^n).$$
Combining this estimate with (4.4) yields (4.3) and completes the proof of (1.4).
\sk

\heading{5. Conclusion}

We have obtained the first two terms in the asymptotic expansions of the moments of $L$
as $\la\to\infty$.
An obvious open question is whether one can obtain a complete asymptotic expansion, or even just
the constant term in (1.3) and the corresponding terms in (1.4).
While our estimates for the contributions to the $O(1)$ term in (1.3) could be improved (for example, by a better choice of the parameter $s$), it is clear that new techniques will be needed to obtain an error term tending to zero.
\sk

\heading{6. Acknowledgment}

The research reported here was supported
by Grant CCF  0917026 from the National Science Foundation.
\sk

\heading{7. References}

\ref C; E. G. Coffman, Jr., T. T. Kadota and L. A. Shepp;
``A Stochastic Model of Fragmentation in Dynamic Storage Allocation'';
SIAM J. Comput.; 14:2 (1985) 416--425.

\refbook N; G. F. Newell;
The $M/M/\infty$ Service System with Ranked Servers in Heavy Traffic;
Springer-Verlag, Berlin, 1984.

\bye